\documentclass[12pt]{amsart}
\usepackage[parfill]{parskip}
\usepackage{graphicx}
\usepackage{wrapfig}
\usepackage{amssymb}
\usepackage{gensymb}
\usepackage{ifxetex}
\ifxetex
  \usepackage{fontspec}
\else
  \usepackage[T1]{fontenc}
  \usepackage[utf8]{inputenc}
  \usepackage{lmodern}
\fi

\usepackage{subfiles}
\usepackage[normalem]{ulem}

\usepackage{amsthm}
\theoremstyle{definition}
\newtheorem{theorem}{{Theorem}}
\newtheorem{lemma}{{Lemma}}

\newtheorem{definition}{{Definition}}
\newtheorem{corollary}{{Corollary}}

\newtheorem{question}{{Question}}

\begin{document}

\title{On Dividing a Rectangle}

\author[Camp Euclid 2016 Group 2]{Robert Dumitru \and Quinn Perian \and Alexander Nealey \and Mohammed Mannan \and Eddie Beck \and Nick Castro \and David Gay \and Dipen Mehta \and Anish Pandya \and Ejay Cho}
\begin{abstract}
In this paper, an interesting geometrical problem has its history reviewed, ranging from the first proof for a partial version by Stanley Rabinowitz, to a new proof of its validity by the authors.
\end{abstract}

\maketitle 

\section{Introduction}

In mathematics, beautiful problems generally have a property in common: that they are simple enough to explain to a person with only an elementary understanding of mathematics, yet in spite of (or perhaps because of) their simplicity, they are very difficult to solve. Our problem is simple enough to explain to a 8 year old, and they would likely come the right answer almost immediately, yet not surprisingly would not know why or how they got that answer, as the proof is not easy to come about. Here is a modified statement of the same problem:
    
\begin{question}
Three brothers inherited a rectangular plot of land from their father, under the condition that they split it so that each brother gets a plot of land with the same measurements as the other two. However, the biggest brother objects to having a rectangular plot of land. Is there any way to split the land such that both the father and the big brother are satisfied?
\end{question}
    
At this point, the authors strongly recommend that the reader try this problem on their own for a few minutes. As is evident to anyone who has tried the problem for a few minutes, it seems impossible to satisfy both the father and the big brother. The difficult part of this problem is not solving it, but in fact proving that there are no solutions.

\section{History}

As far as the authors know, the first proposal of this problem, albeit in a simplified form, was made by Stanley Rabinowitz in 1983, through the problem-solving journal \textit{Crux Mathematicorum} [3]. However, the question asked for a square, not a rectangle. This problem was solved by Samuel J. Maltby in the 1991 volume of \textit{Crux Mathematicorum}[2]. However, it is likely that Rabinowitz had at least a partial solution to the problem, since \textit{Crux} requires submitted problems to have a solution or reasonable evidence that a solution may be arrived at. 
After solving the problem for squares, Maltby went on to prove it for rectangles. As far as the authors know, he was the first person to obtain this result, doing so in 1992[1]. 

There is nothing special about using 3 brothers in the problem. Indeed, we can ask:
\begin{question}
If we replace "Three brothers" with "n brothers", for which n can you satisfy both the father and the big brother?
\end{question}
It is evident that for even n the answer to the question is yes. However, for odd n this is not so clear-cut. It is known that for n greater than 11 the answer is yes. This leaves numbers 3, 5, 7, and 9. 

For n=3, the answer is no, as is proven in this paper and in [1]. For n=5, it is known that the answer is no, although only for squares[4]. It is likely that the answer is no for rectangles as well. Finally, the status of n=7 and n=9 is unknown, as far as the authors know.

Incidentally, the author's attention was drawn to this problem as part of Camp Euclid 2016, under the impression that it was unsolved. The authors came up with a proof of impossibility independently of Maltby. After doing some further research, one of the authors, David Gay, found Maltby's solution. However, this problem gave us much joy in solving it, and the authors felt it would be interesting for the wider mathematical community to know about this problem.

\section{Proof}

Our proof is conceptually quite simple. At heart, it is no more than dividing the problem up into appropriate cases and solving them. However, the reasoning can get quite involved at times, as there are several possible solutions that are valid when rectangles are allowed. To facilitate the discussion of our proof, we introduce some terminology.

\subsection{Definitions}
To facilitate presenting our results, we list some preliminary definitions:

\begin{definition}
By a trisected rectangle we mean a rectangle divided into three non-rectangular regions.
\end{definition}
\begin{definition}
If $S$ is a region of a trisected rectangle, then $\partial S$ denotes the boundary of $S$.
\end{definition}
\begin{definition}
An isometry is a mapping that preserves distances. \\
It can be shown that any isometry is either a translation, rotation, or reflection.
\end{definition}

We will present other definitions as they are needed.

\subsection{Insertion Points}

Perhaps the most important definition that we introduce are insertion points.
\begin{definition}
Suppose that a trisected rectangle has boundary $\Gamma$ and that $S$ and $S'$ are two distinct regions in the trisected rectangle. By an insertion point of $S$ (or $\partial S$), we mean a point that lays on $\Gamma$, and the boundary of $S$. 
\end{definition}

\begin{figure} [h!]
    \centering
    \caption{An example of insertion points.}
	\includegraphics[width=0.50\textwidth]{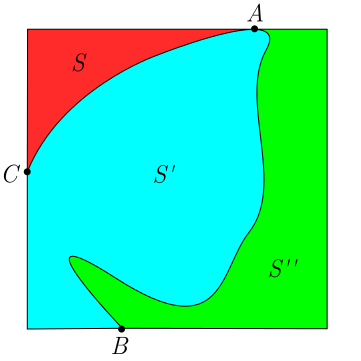}
\end{figure}

In the figure above, $A$ and $C$ are insertion points of $S$; $A$, $B$, and $C$ are insertion points of $S'$; and $A$ and $B$ are insertion points of $S''$.

We approach this problem by classifying trisected squares by the number of insertion points. We then proceed to show that no trisected square with congruent regions falls under any category, which proves that trisected squares with congruent regions do not exist.

\subsection{Types}
We classify regions in a trisected rectangle by the number of insertion points they have:

\begin{itemize}
\item[$\bullet$] Type 0: If a region has no insertion points, then either that region is (a) fully enclosed by the rectangle, or (b) the entire boundary of the rectangle makes up part of the boundary of that region, and the boundary of no other region touches the boundary of the rectangle.
\item[$\bullet$] Type 1: If a region has exactly one insertion point, then that insertion point is either on (a) a corner or (b) an edge.
\item[$\bullet$] Type 2: If a region has exactly two insertion points, then those insertion points are on (a) opposite corners, (b) adjacent corners, (c) opposite edges, (d) adjacent edges, (e) the same edge, (f) a corner and an adjacent edge, or (g) a corner and an opposite edge. Moreover, the two insertion points divide the boundary of the region into two parts. Either one part is fully inside the rectangle and the other part is fully on the boundary of the rectangle, or both parts are fully inside the rectangle.
\item[$\bullet$] Type 3: A region with three or more insertion points. 
\end{itemize}
We say that a region is of type 0, type 1, type 1(a), type 2(b), type 2(c), type 3, etc.

\begin{theorem}
All trisected rectangles with a region of type 2 have a boundary which is split into two parts by the two insertion points, with one part fully inside the rectangle and the other part fully on the boundary of the rectangle.
\end{theorem}
\begin{proof}
If a trisected rectangle has a region of type 2, then the two insertion points divide the boundary of the region into two parts; either one part is fully inside the rectangle and the other part is fully on the boundary of the rectangle, or both parts are fully inside the rectangle. Suppose our rectangle has a region of the latter case. But this region two uniquely determines the other two regions, both of which are of the former case.
\end{proof}

\begin{figure} [h!]
    \centering
    \caption{2 insertion points.}
	\includegraphics[width=0.5\textwidth]{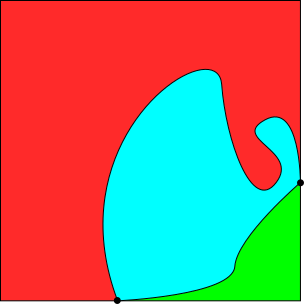}
\end{figure}

\begin{theorem}
 All trisected rectangles with a region of type 3 have a region of type 2.
 \end{theorem}
 
 \begin{proof}
 Call the type 3 region $S$. $S$ has at least 3 consecutive insertion points $A,B$, and $C$ (consecutive in the sense that, along the segment of $\partial S$ containing them, there are no other insertion points). $\partial S$ is divided into three parts by the three insertion points. We look at two parts: the curve connecting $A$ and $B$ and the curve connecting $B$ and $C$. At least one of them is fully interior to the rectangle, which, without loss of generality, we show as the curve connecting $B$ and $C$ (if not, there would either be an insertion point between $B$ and $C$ or $A$ and $B$, or $B$ would not be an insertion point).
 
 We look at how $\partial S$ behaves at the point $C$. First note that the curve connecting $B$ and $C$ divides the rectangle into two parts. $\partial S$ cannot go into the part not containing $A$, otherwise it cannot connect to $A$ without intersecting itself. Furthermore, if $\partial S$ connects to $A$, it must be the case that another region is created in addition to $S$ and the region whose boundary is made up of the curve interior to the rectangle connecting $B$ and $C$ and the part of the boundary of the rectangle connecting $B$ and $C$. However, the latter region is of type 2. Therefore, all trisected rectangles with a region of type 3 must have a region of type 2.
 \end{proof}
 
 \begin{figure} [h!]
    \centering
    \caption{3 insertion points.}
	\includegraphics[width=0.50\textwidth]{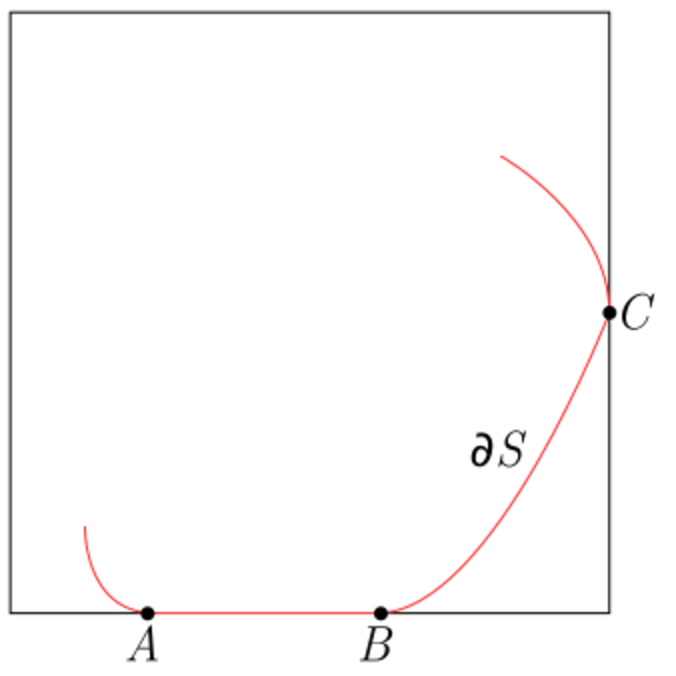}
\end{figure}

\begin{theorem}
All trisected rectangles with a region of type 0 have a region of type 1 or 2.
\end{theorem}
\begin{proof}
Suppose that a trisected rectangle has a region of type 0(b). Then it has a region $S$ that has opposite corners of the rectangle, i.e., it has two points of maximum possible distance in a rectangle. Since the other two regions have no points on the boundary of the rectangle, they cannot have this property. Hence they cannot be congruent to $S$, so there are no trisected rectangles with a region of type 0(b). 

Now, suppose that a rectangle has a region of type 0(a). Then it has a region $S$ that is fully enclosed by the rectangle. We look at the two other regions S' and S". It cannot be that the boundary of the rectangle is fully on $\partial$ S' and doesn't touch $\partial S''$, or vice versa (since then the rectangle would have a region of type 0(b)). Hence there is a point on the boundary of the rectangle where $\partial S'$ and $\partial S''$ coincide. This is an insertion point. Hence the rectangle contains a region with multiple insertion points. So it either has 1, 2, or more than 2 insertion points. But if it has a region with more than 2 insertion points, by Theorem 3, it has a region with 2 insertion points. Hence the trisected rectangle has a region of type 1 or 2.
\end{proof}
\begin{corollary}
All trisected rectangles have a region of type 1 or 2.
\end{corollary}
\begin{proof}
Follows from Theorems 1-3.
\end{proof}

\begin{theorem}
No trisected rectangle with congruent regions contains a region of type 1(a).
\end{theorem}
\begin{proof}
Suppose that there is a trisected rectangle with a region of type 1(a). Then it has a region $S$ whose boundary touches the boundary of the rectangle at only one point, the corner $A$. The other two regions, $S'$ and $S''$, cannot contain opposite corners (otherwise they would contain two points whose distance from each other is larger than the distance between any two points of $S$). Hence $S'$ contains $A$ and $B$ and $S''$ contains $C$ and $D$, or $S'$ contains $A$ and $D$ and $S''$ contains $B$ and $C$. The two cases are equivalent, being reflections of each other, so without loss of generality we consider the former case.

By definition, $S'$ must contain $S$. Since $S$ must map to $S'$, $S$ must contain another region, $S_1$. We can apply the same logic ad infinitum, thereby giving us a infinity of regions. This contradicts the fact that a trisected rectangle contains 3 regions. Therefore, a trisected rectangle cannot contain a region of type 1(a).
\begin{figure} [h!]
    \centering
    \caption{A region of type 1(a).}
	\includegraphics[width=0.50\textwidth]{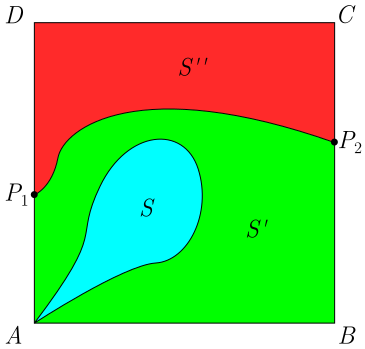}
\end{figure}
\end{proof}
\begin{theorem}
No trisected rectangle with congruent regions has a region of type 1(b).
\end{theorem}
\begin{proof}
Suppose that there exists a trisected rectangle with a region of type 1(b). Then it contains a region $S$ such that $\partial S$ touches the boundary of the rectangle at one point on the edge, $P$. The other two regions, $S'$ and $S''$, cannot contain opposite corners (otherwise they would contain two points whose distance from each other is larger than the distance between any two points of $S$). Hence \textbf{(I)} $S'$ contains $A$ and $B$ and $S''$ contains $C$ and $D$, or \textbf{(II)} $S'$ contains $A$ and $D$ and $S''$ contains $B$ and $C$.

\textbf{(I)}: $\partial S'$ contains the entire edge $AB$ and $\partial S''$ contains the entire edge $CD$ (otherwise there would be more than 3 regions). Furthermore, $\partial S'$ contains a portion of $AD$ adjacent to $A$ (otherwise $S''$ contains $A$ as well). Similarly, $\partial S'$ contains a portion of $BC$ adjacent to $B$, $\partial S''$ contains a portion of $AD$ adjacent to $D$, and $\partial S''$ contains a portion of $BC$ adjacent to $C$. Then $\partial S'$ has an insertion point $P_1$ on $AD$ and $P_2$ on $BC$. $P_1$ and $P_2$ divide $\partial S'$  into two parts. Consider the part not containing $AB$. If it does not intersect $\partial S$, then $\partial S'$ has a point where it intersects itself, $P$. $\partial S''$ does not have this property, so $S'$ and $S''$ cannot be congruent. Otherwise, it intersects $\partial S$. But this cannot be done without having more than three regions in the rectangle.

\textbf{(II)}: $\partial S'$ contains the entire edge $AD$ and $\partial S''$ contains the entire edge $BC$ (otherwise there would be more than 3 regions). Furthermore, $\partial S'$ contains a portion of $AB$ adjacent to $A$ (otherwise $S''$ contains $A$ as well). Similarly, $\partial S'$ contains a portion of $CD$ adjacent to $D$, $\partial S''$ contains a portion of $AB$ adjacent to $B$, and $\partial S''$ contains a portion of $CD$ adjacent to $C$. Then $\partial S'$ has an insertion point $P_1$ on $AB$ and $P_2$ on $CD$. $P_1$ and $P_2$ divide $\partial S'$ into two parts. Consider the part $Z$ not containing $AD$. If $P_1\neq P$, then one of $\partial S'$ or $\partial S''$, intersects itself at $P$, while the other has no point where it intersects itself. Hence $P_1=P$.

If $Z$ only intersects $\partial S$ at $P_1$, then either $\partial S'$ or $\partial S''$ intersects itself at $P$ while the other has no point where it intersects itself, so $Z$ intersects $\partial S$ at $P_3\neq P_1$. The part of $Z$ between $P_1$ and $P_3$ must be the same as the part $\partial S$ between $P_1$ and $P_3$ on the side of $AD$. This uniquely determines the region $S'$, which in turn uniquely determines the region $S''$.

Since $S'$ and $S''$ are congruent, there is a isometry mapping $S'$ to $S''$. Under this mapping, $AD$ maps to a line segment on $\partial S''$. But this cannot be $AP_1$ or $DP_2$, nor can it be any part of $Z$. Hence $AD$ maps to $BC$. Furthermore, the map cannot send $A$ to $C$ and $D$ to $B$. Therefore, $S''$  is the reflection of $S'$ over the vertical line splitting the rectangle in half. Thus, $P_2$ is the midpoint of $DC$, $ P_1$ is the midpoint of $AB$, and the part of $Z$ between $P_2$ and $P_3$ is a vertical line segment. Moreover, $S$ has reflection symmetry about the line segment $P_{1}P_{2}$. $S'$ and $S''$ do not have this property(they have no reflection symmetry at all), so they cannot be congruent to $S$. Thus, there exist no trisected rectangles with congruent regions with a region of type $1(b)$. 
\begin{figure} [h!]
    \centering
	\includegraphics[width=0.50\textwidth]{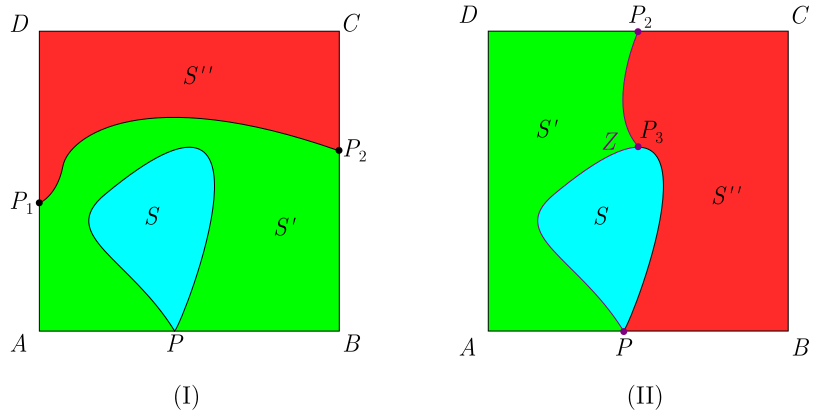}
\end{figure}
\end{proof}
\begin{corollary}
There are no trisected rectangles with congruent regions with a region of type $1$.
\begin{proof}
Follows from Theorems $4$ and $5$. 
\end{proof}
\end{corollary}
\begin{theorem}
There are no trisected squares with three congruent regions that contain a region of type 2(a) (a region with exactly 2 insertion points such that the insertion points are on opposite corners).
\end{theorem}
\begin{proof}
Suppose that there exists a trisected rectangle with a region of type $2(a)$. Then it contains a region $S$ such that $\partial S$ is formed by joining a curve connecting opposite corners $A$ and $C$ interior to the rectangle with sides $AB$ and $BC$. Since $S$ has two opposite corners, it has two points with the maximum possible distance between two points in a rectangle. The other two regions, $S'$ and $S''$, must also have this property. The point $B$ cannot belong to $S'$ or $S''$. Hence $S'$ and $S''$ must have $A$ and $C$.Since $S$, $S'$, and $S''$ are congruent, there exists a rigid motion mapping $S$ to $S''$. But such rigid motion fixes distances, so it either maps $A$ to 
$A$ and $C$ to $C$, or $A$ to $C$ and $C$ to $A$. Hence $S''$ is either \textbf{(I)} the reflection of $S$ over $AC$ or \textbf{(II)} the rotation of $S$ about the center of the rectangle through an angle of $\pi$.The boundary of $S'$ is the joining of the parts boundary of $S$ and $S'$ that are interior to the rectangle. 

\textbf{(I)} $S'$ has reflection symmetry about $AC$. But $S$ does not have reflection symmetry about its line segment of maximum distance, $AC$. Hence $S$ and $S'$ cannot be congruent. 

\textbf{(II)} $S'$ has rotational symmetry through an angle of $\pi$ about the center of its line segment of maximum distance. But $S$ does not have this property. Since both cases are impossible, a trisected rectangle with a region of type $2(a)$ cannot exist.
\end{proof}
\begin{figure} [h!]
    \centering
	\includegraphics[width=0.60\textwidth]{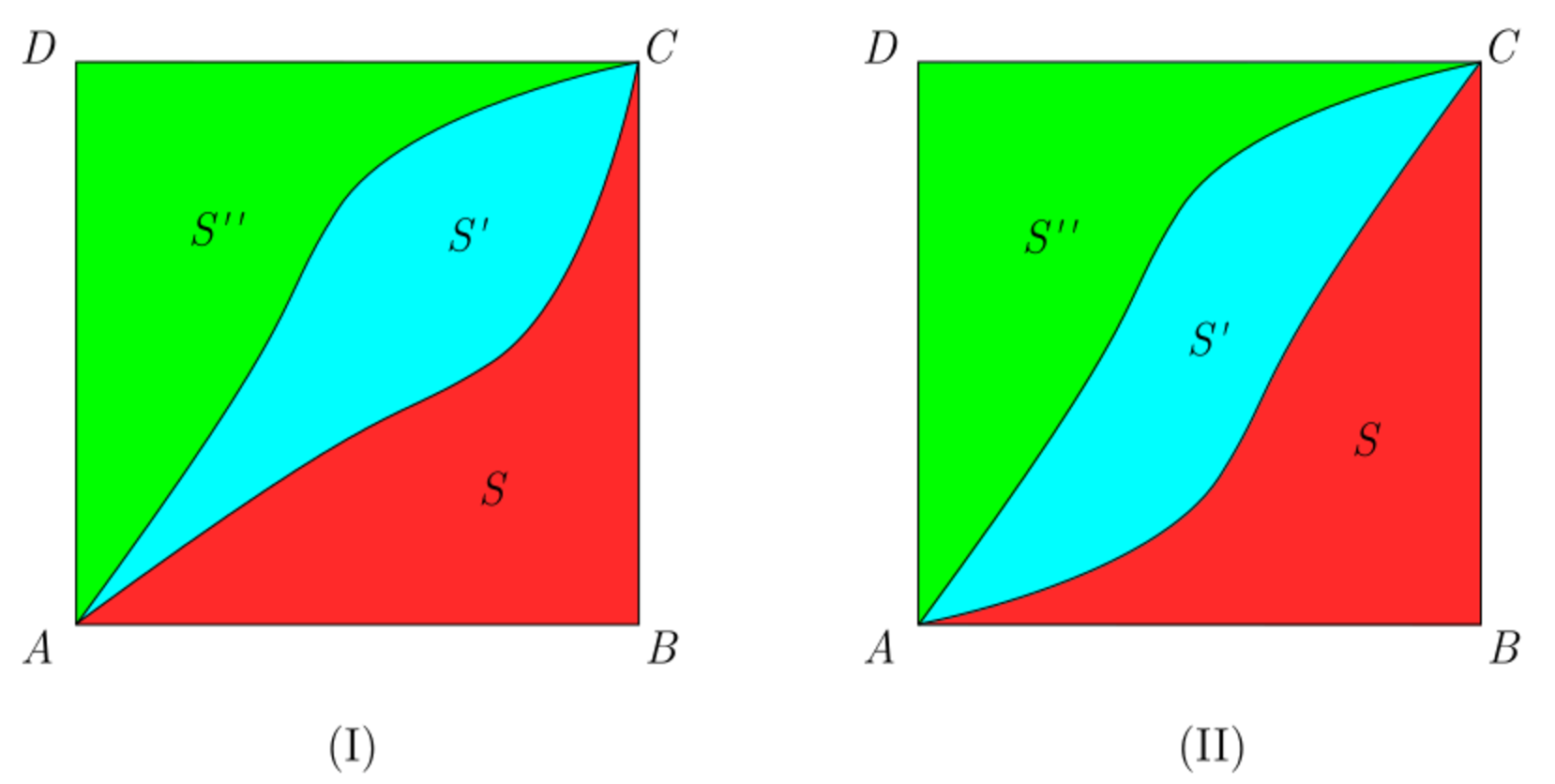}
\end{figure}
\begin{theorem}
There are no trisected squares with three congruent regions that contain a region of type 2(b) (a region with exactly 2 insertion points such that the insertion points are on adjacent corners).
\end{theorem}
\begin{proof}
Suppose that a trisected rectangle of type 2(b) exists. Then for one shape to not have two opposite corners, there are four possible configurations:

\textbf{(I)}: Line $Y$ goes from corner $D$ to edge $AB$. If this is the case, then shape $S$ needs to have more straight lines that the minimum because it has at least straight line (line $CD$), and shape $S'$ and $S''$ each have at least two straight lines(lines $AC + F$ and lines $BD + E$ respectively) while not having any of these straight lines touch other shapes, as otherwise those shapes will still have more straight lines in comparison. This is impossible as to create a straight line in shape $R$, $t$ must go somewhere on curve $X$, but curve $X$ touches a different shape at all points. 

\begin{figure} [h!]
    \centering
	\includegraphics[width=0.60\textwidth]{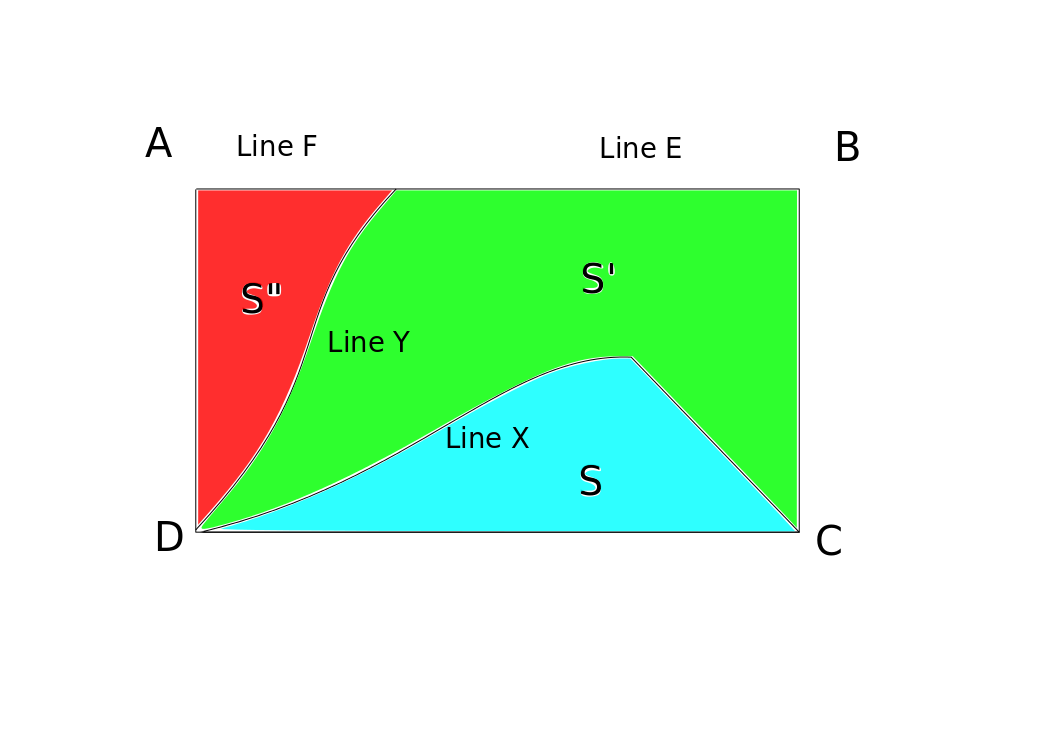}
\end{figure}

\textbf{(II)}: Line $Y$ goes from corner $D$ to edge $AB$. This case is the same as case (I) except that the diagram is reflected over the x-axis. Therefore, the same proof applies here.

\textbf{(III)}: Line $Y$ goes from edge $AC$ to edge $CD$. In this case if line $X$ has $X$ straight lines and line $Y$ has $Y$ straight lines, shape $R$ has 1(line $CD$) + $X$ straight lines, shape $R'$ has 2(lines $G$ and $H$) + $X$ + $Y$ straight lines, and shape $R''$ has 3(lines $E$, $F$, and $AB$) + $Y$ straight lines. Therefore $1+ X = 2 + X + Y = 3 + Y$, which is impossible as $2 + X + Y > 1 + X$.

\begin{figure} [h!]
    \centering
	\includegraphics[width=0.60\textwidth]{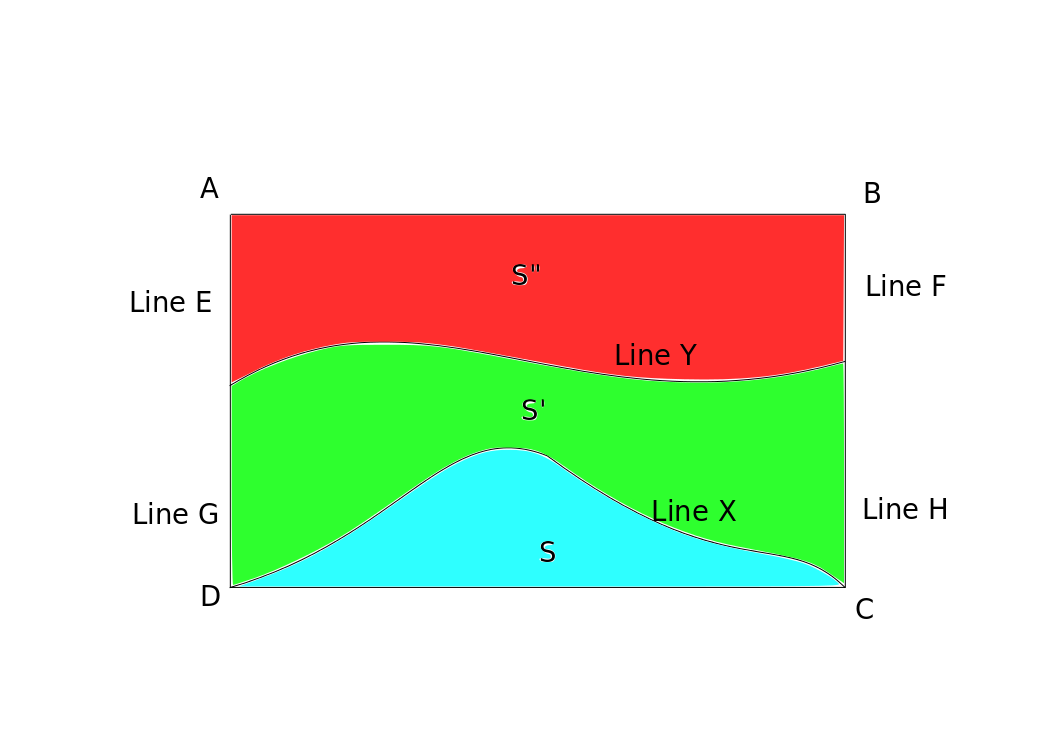}
\end{figure}

\textbf{(IV)}: Line $Y$ goes from edge $AB$ to curve C.  Since S is congruent to S', DC must map to PQ, BC, CP and BQ. 

\textbf{(i)}: DC maps to BC. since $BC=DC$,  $ABCD$ is a square. In this case all shapes other than shape $S$ have a right angle that is made up of the boundary of the square. The only place where the right angle can map to is somewhere on line $X$, as the other two angles each are only part of a $90^{\circ}$ angle. If this were the case then the side opposite line $X$, which is on the boundary of the rectangle would have to map to the longest curve(non-straight line segment)/“hypotenuse” of one of the other shapes, but this is impossible as that line is longer than one side of the square. This line has to be the longest line of the shape as if it were a right triangle that line would be that hypotenuse, and because a straight line is the shortest distance between two point, you can't make that line shorter. Therefore, DC cannot map to BC.

\textbf{(ii)}: DC maps to CP. Then, by symmetry, DC also maps to DP. Since $DC=CP=DP$, S is an equilateral triangle. Since S is congruent to S' and S", S' and S" are also equilateral. However, both S' and S" contain $90^{\circ}$angles. Therefore, DC cannot map to CP.

\textbf{(iii)}: DC maps to PQ. Then, since S is congruent to S',  either $\angle{DCP}$ or $\angle{CDP}$ is equal to $90^{\circ}$. This  is absurd if we intend to keep S a region of type 2(b). Therefore, DC cannot map to PQ.

\textbf{(iv)}: DC maps to BQ. However, that implies that $BQ=AB$, which is absurd.
\begin{figure} [h!]
    \centering
	\includegraphics[width=0.60\textwidth]{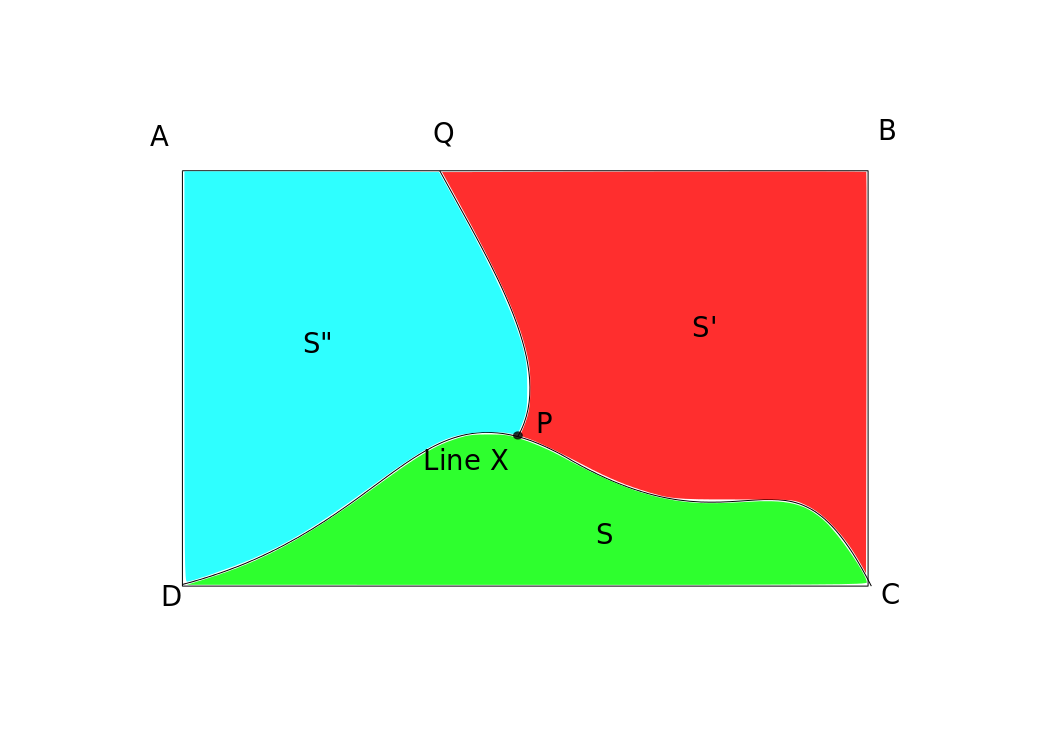}
\end{figure}

Therefore, as all cases are impossible, a trisected rectangle with congruent regions cannot have a region of 2(b).
\end{proof}

\begin{theorem}
There are no trisected squares with congruent regions with a region of type 2(c)(a region with exactly 2 insertion points such that the insertion points are on opposite edges).
\end{theorem}
\begin{proof}
 For this case we are assuming we are working with squares.
 We will consider a rectangle with a 2(c) region and exhaustively consider every possible way of splitting the remainder of the rectangle into two regions. To streamline this proof, we propose the following definitions for this theorem only:
 
 \begin{definition}
 By a 2(c) boundary, we mean the part of the boundary that is not part of the boundary of the rectangle.
 \end{definition}
\begin{definition}
By a partition, we mean the line that separates the remainder of the rectangle into two regions.
\end{definition}
\begin{definition}
By an adjacent edge, we mean an edge that is connected to the 2(c) boundary.
\end{definition}
\begin{definition}
By a opposite edge, we mean an edge that is not connected to the 2(c) boundary.
\end{definition}

We consider all possible partitions of the remainder of the rectangle into two regions. These fall under three broad categories: \textbf{(I)}, where both endpoints on the 2(c) boundary; \textbf{(II)}, where only one endpoint is on the 2(c) boundary; and \textbf{(III)}, where no endpoints are on the 2(c) boundary. With two endpoints, this list is exhaustive.

\textbf{(I)}: We have three cases to consider: both endpoints are on the interior of the 2(c) boundary, which is impossible because of lemma 2; only one endpoint is on the interior, which is impossible by by corollary 2; and neither endpoint is on the interior of the 2(c) boundary, which is impossible because of lemma 1. 

\begin{lemma}
A region of type 2(c) is not possible when the remainder of the rectangle is separated by a partition where neither endpoint is on the interior of the 2(c) boundary. 
\end{lemma}
\begin{proof}
Observe that EB and FA are parallel, distance EF apart, and, if extended to full lines, would completely contain region BEFA. However, for there to be two segments with the same properties in the boundary of the center region, we would need to have the segment that FA maps to parallel to FD because the central region extends a width of EF from EC to FD. However, this would clearly make one of the outer regions a rectangle. Thus, this case is impossible.\\
\begin{center}
\includegraphics[width=\textwidth]{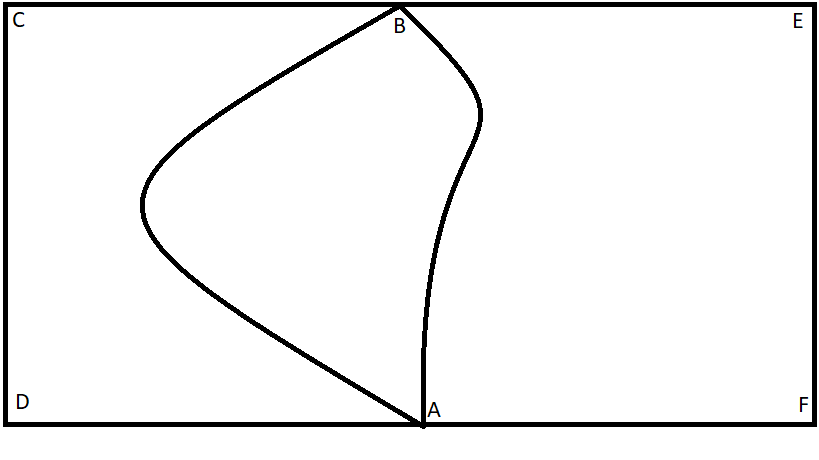}
\end{center}
\end{proof}
\begin{lemma}
A region of type 2(c) is not possible when the remainder of the rectangle is separated by a partition where both end points are on the interior of the 2(c) boundary. 
\end{lemma}
\begin{proof}
We will assume that these two endpoints are distinct. If this is not the case, then region $S''$ clearly can't map to region $S$. We have two ways to map S to S": \textbf{(i)}AP maps to QC; or \textbf{(ii)}AP maps to DP. 

\textbf{(i)}: It immediately follows that S' has $180^{\circ}$ symmetry. However, for S and S" to have the same property, the 2(c) boundary must be a straight line, which makes both S and S'' a rectangle. However, this clearly violates the conditions of the problem. Therefore, AP cannot map to QC.

\textbf{(ii)}: It immediately follows that S' has reflectional symmetry. Since S and S" are isometries of S', they must also have reflectional symmetry. Therefore, $S\cup S''$ must have two axes of symmetry, each perpendicular to each other. From this we can infer that S' also has two perpendicular axes of symmetry. However, S and S" must also the same property, which forces them to be rectangles. This is forbidden by the problem statement. Therefore, including the conclusion from \textbf{(i)}, we cannot have a region of type 2(c) when the remainder of the rectangle is separated by a partition where both end points are on the interior of the 2(c) boundary. 
\begin{figure}[h!]
    \centering
	\includegraphics[width=0.60\textwidth]{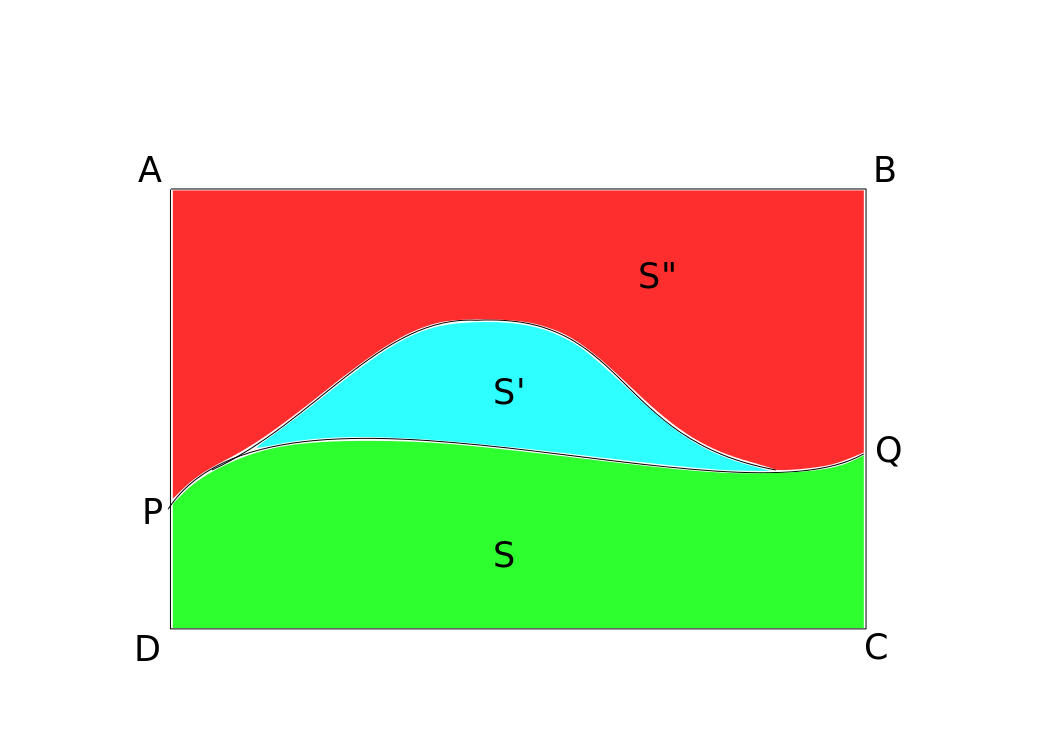}
\end{figure}
\end{proof}
\textbf{(II)}: There are six cases to consider: (i) Three where the endpoint is on the interior of the 2(c) boundary; (ii) three where it is on the exterior. In each of those three cases the endpoint is on the adjacent edge, the corner or the opposite edge. 

\textbf{(i)}: In the case where one endpoint is on the interior of the 2(c) boundary, the other endpoint cannot be on the adjacent edge by theorem 10 and cannot go to the corner by both theorems 11 and 12. Therefore, we need only prove the case in which the other endpoint is on the exterior.

\begin{lemma}
A region of type 2(c) is not possible when the remainder of the rectangle is separated by a partition where the partition has one endpoint on the 2(c) boundary and the other on the opposite edge. 
\end{lemma}
\begin{proof}
There are 6 cases: (I) EF and PQ are both straight, (II) EF is straight but not PQ, (III) PQ  is straight but not EF , (IV) Both EF and PQ are not straight , (V) EQ and PQ are straight but not QF , and (VI) QF is straight but PQ and EQ are not.

\textbf{(I):} Assume WLOG that right angle $\angle{CDE}$ maps to right angle $\angle{EAP}$. Then, either EF is parallel to AB, or PQ is parallel to AD. If EF is parallel to CD, then CDEF is a rectangle, which is forbidden. If PQ is parallel to AD, then $[APQE]\neq[BFQP]$. The only time $[APQE]$ = $[BFQP]$ is when all 3 regions are rectangles.

\textbf{(II):} $CDEF$ has 4 straight lines as boundaries. Therefore, $APQE$ and $BFQP$ will also have 4 straight lines as boundaries. WLOG, let us examine $APQE$. QE, AE, and AP are all straight lines. Therefore, PQ must also be a straight line to complete the boundary. However, we expressly assumed that PQ is not straight.

\textbf{(III):} As in case (I), assume WLOG that right angle $\angle{CDE}$ maps tho right angle $\angle{EAP}$. Then, DC must map AP. However, that is impossible as $AP<AB=DC$.

\textbf{(IV):}CD cannot map to AP, for $CD>AP$.  Therefore, AE must map to CD. However, CF must map to EQ, which is a contradiction.

\textbf{(V):}EF is evidently not straight, as QF is not straight. However, that implies that EF must map to one of $APQE$'s straight sides, which a contradiction.

\textbf{(VI):} For the same reason as in case (IV), CD must map to AE. However, CF must map to EQ, which is a contradiction.
\newpage
\begin{figure}[h!]
    \centering
	\includegraphics[width=0.60\textwidth]{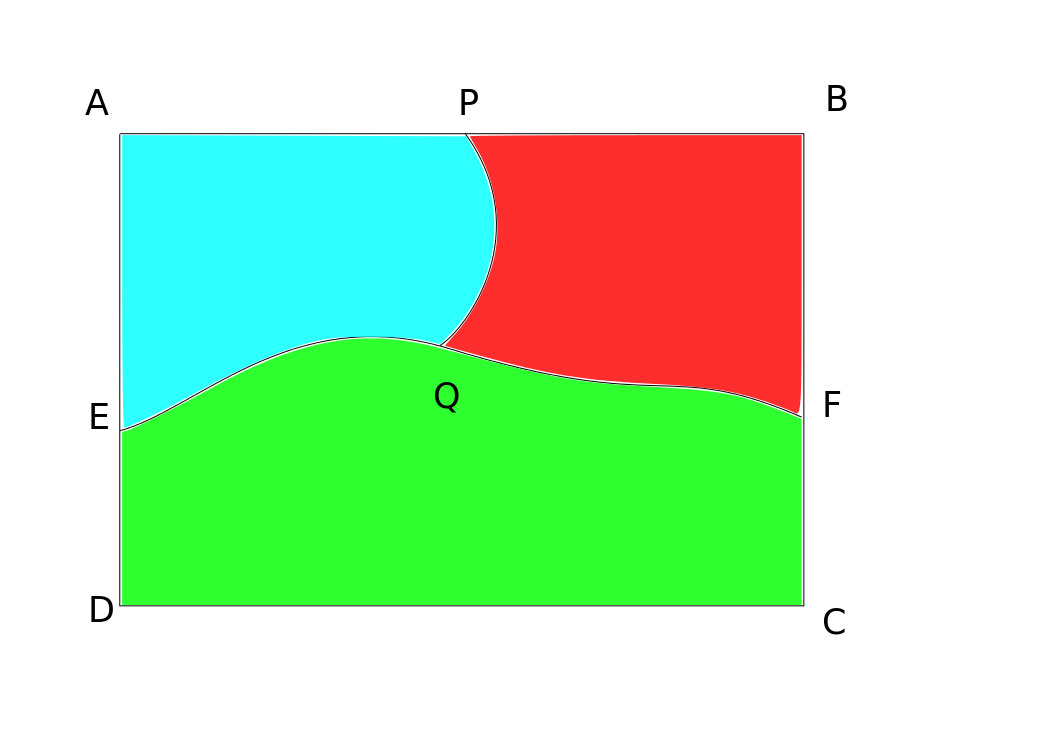}
\end{figure}
\end{proof}
\textbf{(ii)}: The partition cannot go from the exterior of 2(c) boundary to the corner by theorem 12 and can't go to the opposite edge by the same logic as lemma 3. 

\begin{lemma}
A region of type 2(c) is not possible when the remainder of the rectangle is separated by a partition where the partition has one endpoint on the exterior of the 2(c) boundary and the other on the other adjacent edge. 
\end{lemma}
\begin{proof}
This case has two regions of 2(c) (S and S'') and one region of type 3 (S'). S and S'' have two straight edges that are parallel to each other and distance of the height of the rectangle away from each other.

If the rectangle is a square, the straight edge on S' can't be on the parallel edges as then it would no longer be a type 3 region. And the other two edges can't be the parallel edges as their distance from each other is strictly less than the width of the square.

If the rectangle is not a square, we know that S' spans the height of the rectangle while the sum of the lengths of the opposite edges have to equal the width. The forces S' to be a type 2 region specifically one disproved by lemma 1.

Therefore, a region of type 2(c) is not possible when the remainder of the rectangle is separated by a partition where the partition has one endpoint on the exterior of the 2(c) boundary and the other on the other adjacent edge.  
\begin{figure}[h!]
    \centering
	\includegraphics[width=0.50\textwidth]{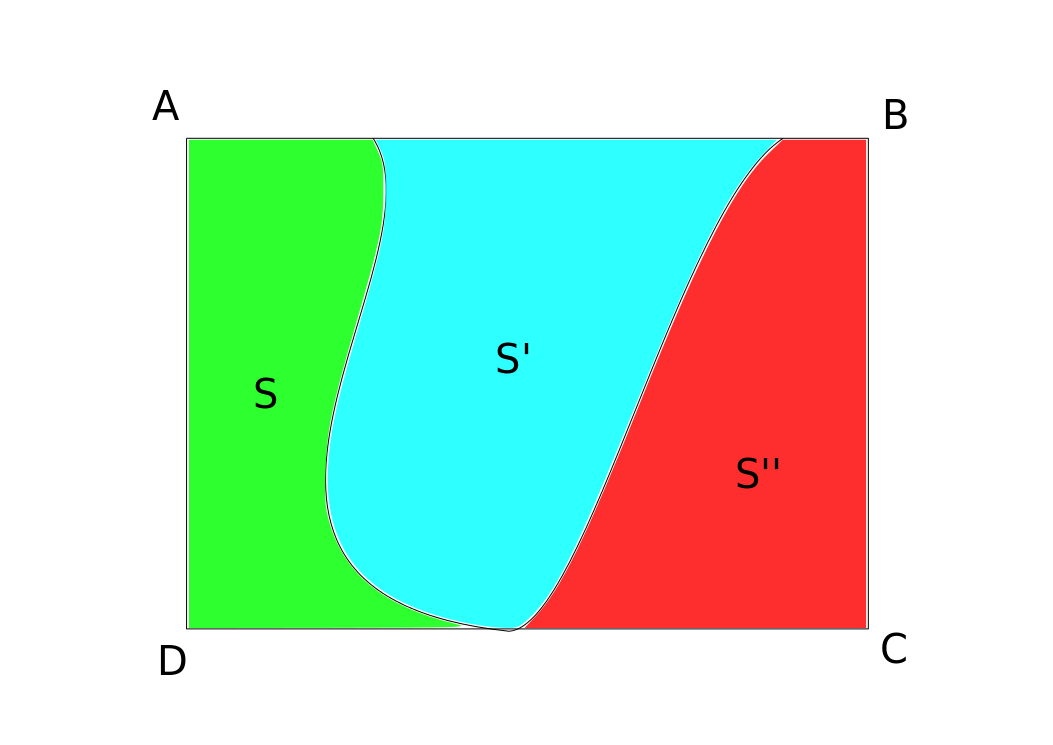}
\end{figure}
\end{proof}

\textbf{(III)}: This can be considered in three cases: (i), where at least one endpoint of the partition is on an adjacent edge; (ii), where both endpoints on the opposite edge; and (iii), no endpoints are on any edge.  

\textbf{(i)}:There are five cases for us to consider here: where the other endpoint is on the same adjacent edge (impossibly by theorem 10); on the corner of the same adjacent edge (impossible by theorem 11); on the corner of the other adjacent edge (impossible by theorem 12); on the interior of the other adjacent edge (lemma 5); or on the interior of the opposite edge (lemma 6).  

\begin{lemma}
A region of type 2(c) is not possible when the remainder of the rectangle is separated by a partition where the partition has one endpoint on each adjacent edge. 
\end{lemma}
\begin{proof}
Each of the outer regions must have three straight edges that meet at right angles, but the middle region only starts with two. Therefore, one of $M$ and $N$ must be straight and perpendicular to $AB$ and $CD$, because the middle region maps to the outer ones. However, this forces one of the outer regions have 4 straight sides and 4 right angles, which is the definition of a rectangle.
\begin{figure}[h!]
    \centering
	\includegraphics[width=0.50\textwidth]{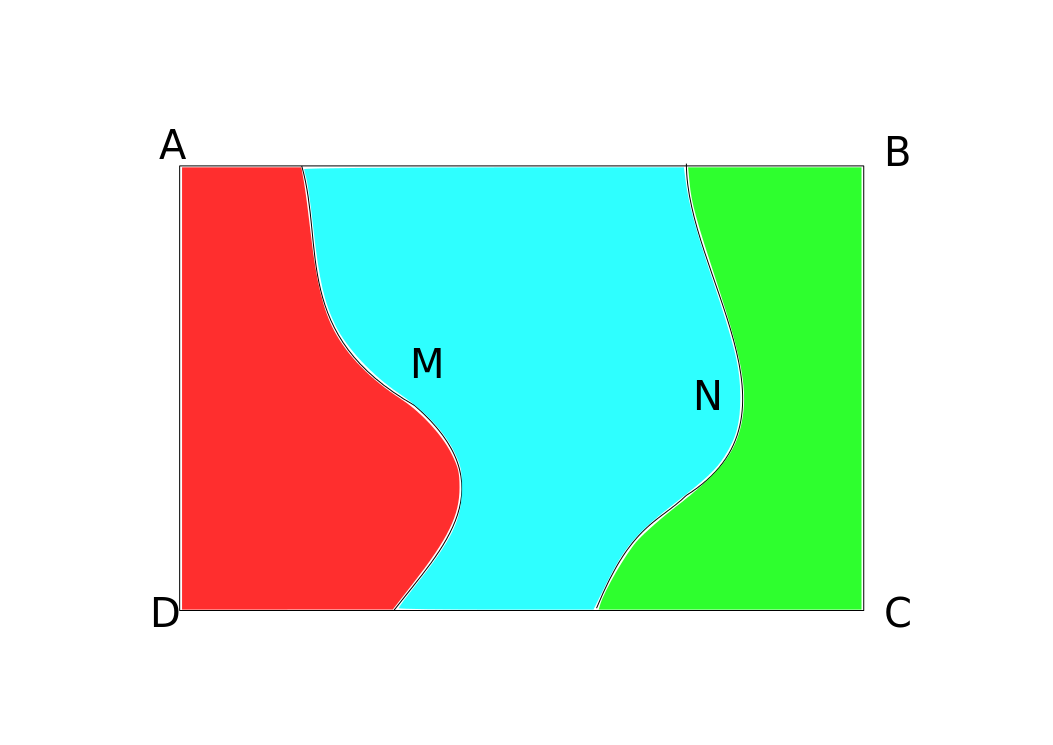}
\end{figure}
\end{proof}
\newpage
\begin{lemma}
A region of type 2(c) is not possible when the remainder of the rectangle is separated by a partition where the partition has one endpoint on an opposite edge, and one on an adjacent edge. 
\end{lemma}
\begin{proof}
In this case, $S''$ starts with only two straight lines, and the remaining line touches $S'$, which has 3 straight lines. Therefore, if  $S''$ has another straight line, $S'$ has another straight line as well, making the $S''$ always have less straight lines. Therefore, $S''$ is not congruent to $S'$.
\begin{figure}[h!]
    \centering
	\includegraphics[width=0.60\textwidth]{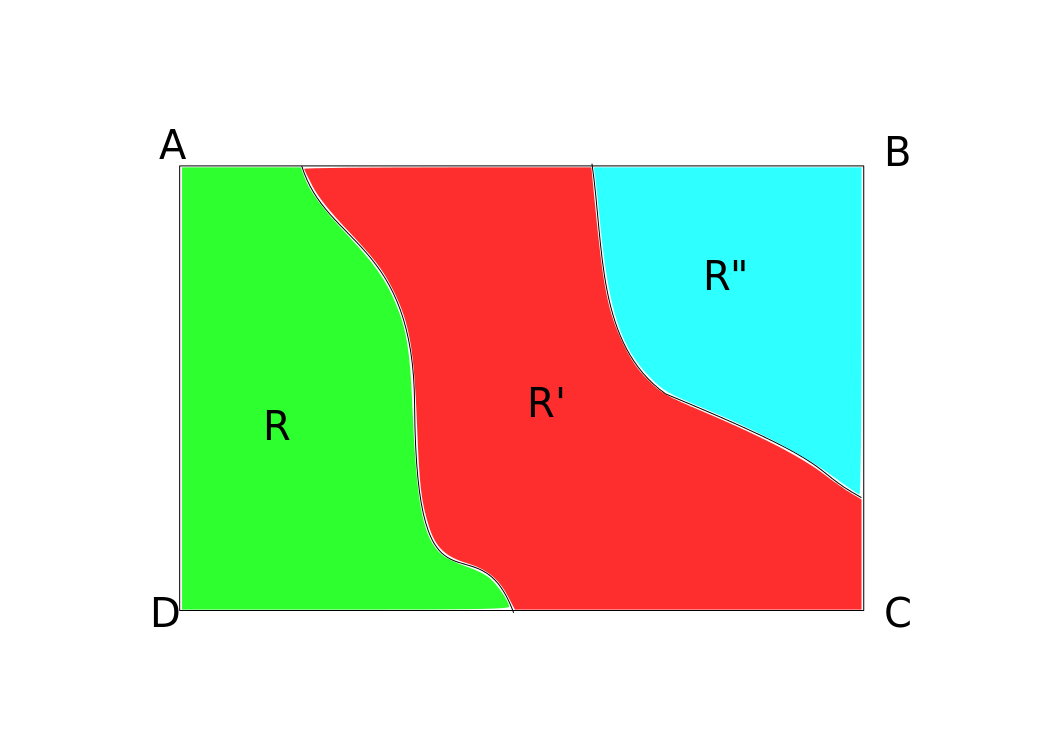}
\end{figure}
\end{proof}

\textbf{(ii)}:This can be considered in three cases: Both endpoints are on the corners (which is not possible by theorem 6), both endpoints are on the interior of the opposite edge (which is not possible by theorem 6) and one endpoint is on the interior while the other is on the corner (which is not possible by theorem 11).  

\textbf{(iii)}:
\begin{lemma}
The endpoints of the partition must be on an edge. 
\end{lemma}
\begin{proof}
Assume that no endpoints of the partition are on any edge. By the hypothesis of the problem, $R'$ is congruent to $R''$. However, since $R'$ contains a region($R''$), $R'$ must have a region inside it. However, there are four regions in this rectangle, while the problem asks for three, so by reductio ad absurdum, the lemma is proved.
\begin{figure}[h!]
    \centering
	\includegraphics[width=0.60\textwidth]{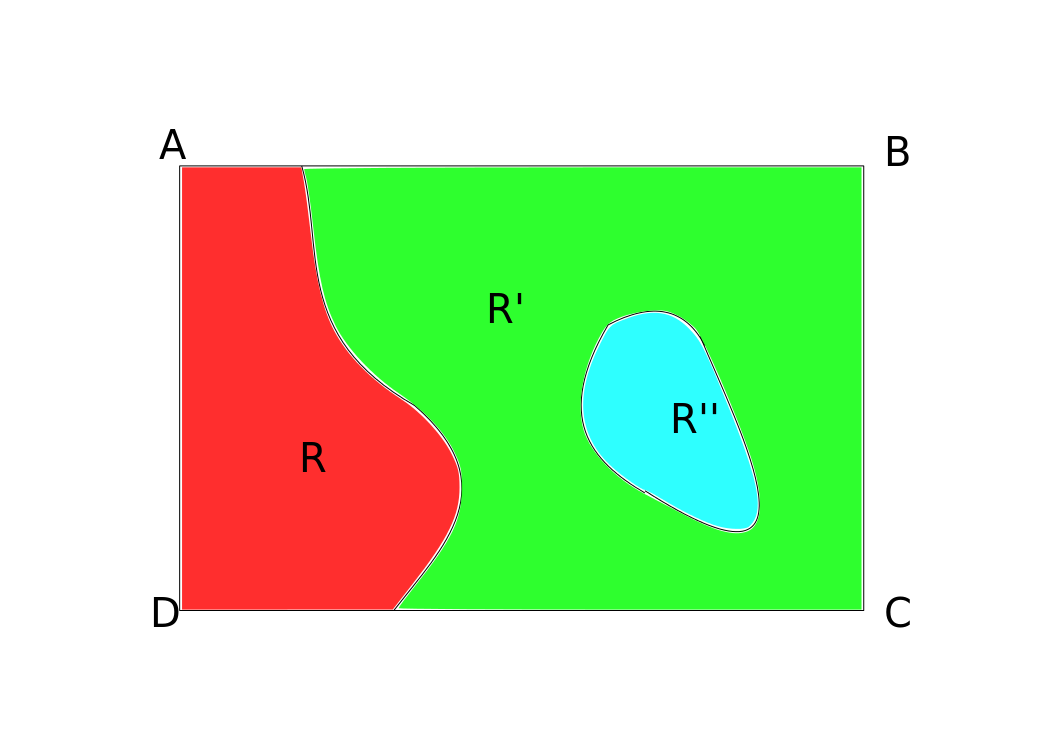}
\end{figure}
\end{proof}
\end{proof}
\begin{theorem}
There are no trisected squares with congruent regions with a region of type 2(d) (a region with exactly 2 insertion points such that the insertion points are on adjacent edges).
\end{theorem}
\begin{proof}
Suppose that a trisected rectangle has a region of type 2(d). In all images, the region of type 2(d) will be on the bottom right corner. Then there are three cases of such a trisected rectangle: (I) one region contains opposite corners; (II) the second line does not start or end on a corner; and (III) the second line starts or/and ends on a corner. 

\textbf{(I)}: Case I is impossible because the region with opposite corners contains a main diagonal of the rectangle, which the other two regions do not have. Therefore, the three regions cannot be congruent. 

\textbf{(II)}: For the three regions to be congruent, all three must have the same number of straight edges. However, in all the cases, $S'$ has the most straight edges. This is because $S'$ always shares an edge with another region. If we make the internal edge straight, it will give $S$ or $S''$ another straight edge, but at the same time, it will also give $S'$ another straight edge. Therefore, the three regions cannot be congruent. Note that this argument does not hold for the third example, because it has a place that only the red and green areas share an internal line.  It is included below.  

\begin{figure}[h!]
    \centering
	\includegraphics[width=0.60\textwidth]{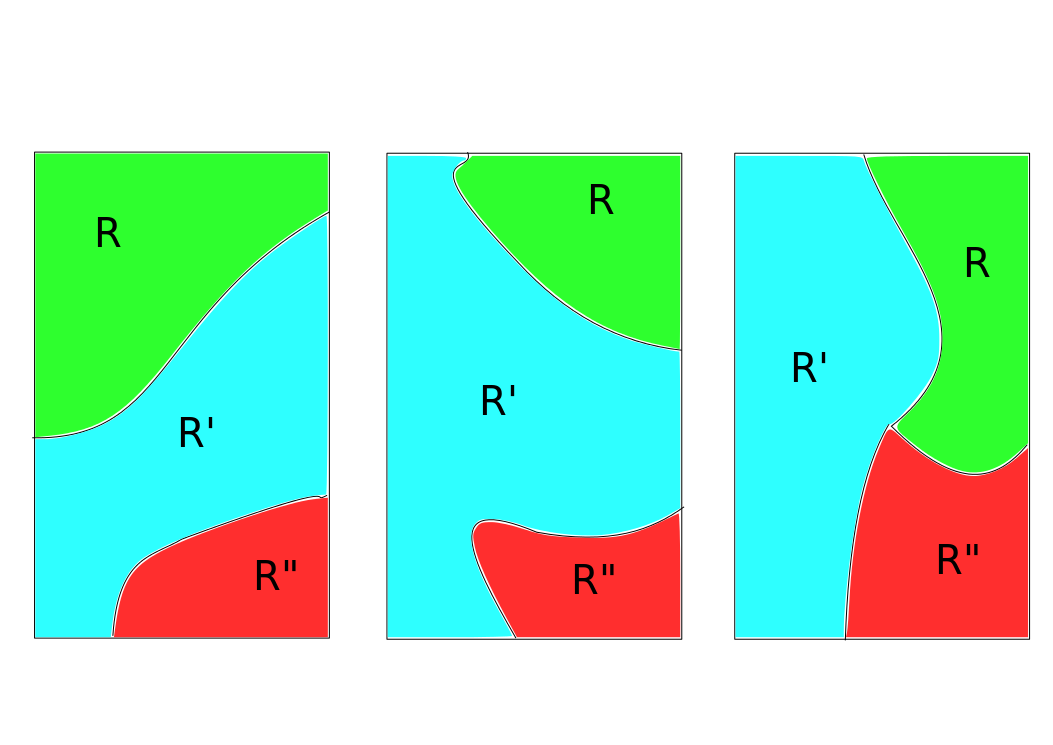}
\end{figure}

\textbf{(III)}: For the three regions to be congruent, all three must have the same number of straight edges. However, in the first seven cases, $S'$ has the most straight edges because, by adding a straight line to $S$ or $S''$, we must also add a straight line to $S'$. Therefore, the three regions cannot be congruent in those seven cases. 

\begin{figure}[h!]
    \centering
	\includegraphics[width=0.80\textwidth]{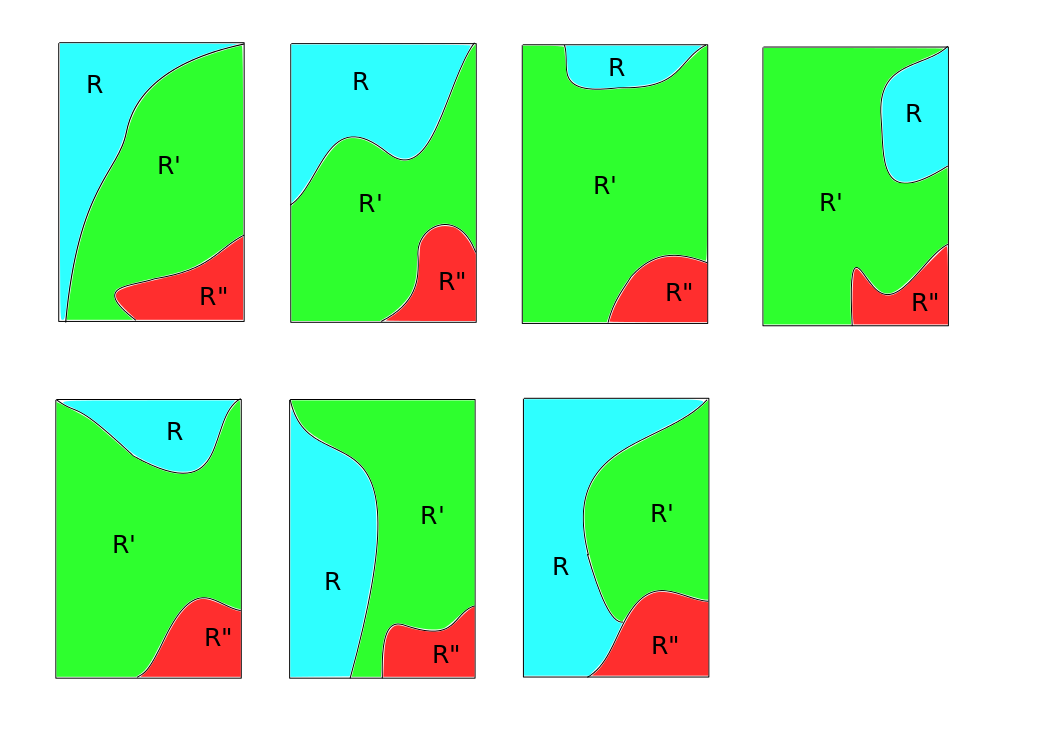}
\end{figure}

Note that we left two cases. These two cases are special in that all three regions can have the same number of straight edges. However, Special Case A is still impossible because it contains a region of type 2(a). This is impossible according to Theorem 6. 

This leaves us with the cases below. These cases all fall under a category we have already proved, cases 1, 4, 5, and 6 are each of type 2(g), while cases 2, 3, and 7 are also of type 2(c). Therefore, as all cases of this type are impossible, there can be no shape of type 2(d) in a trisected square.

\end{proof}

\newpage
\begin{figure} [h!]
    \centering
	\includegraphics[width=1.05\textwidth]{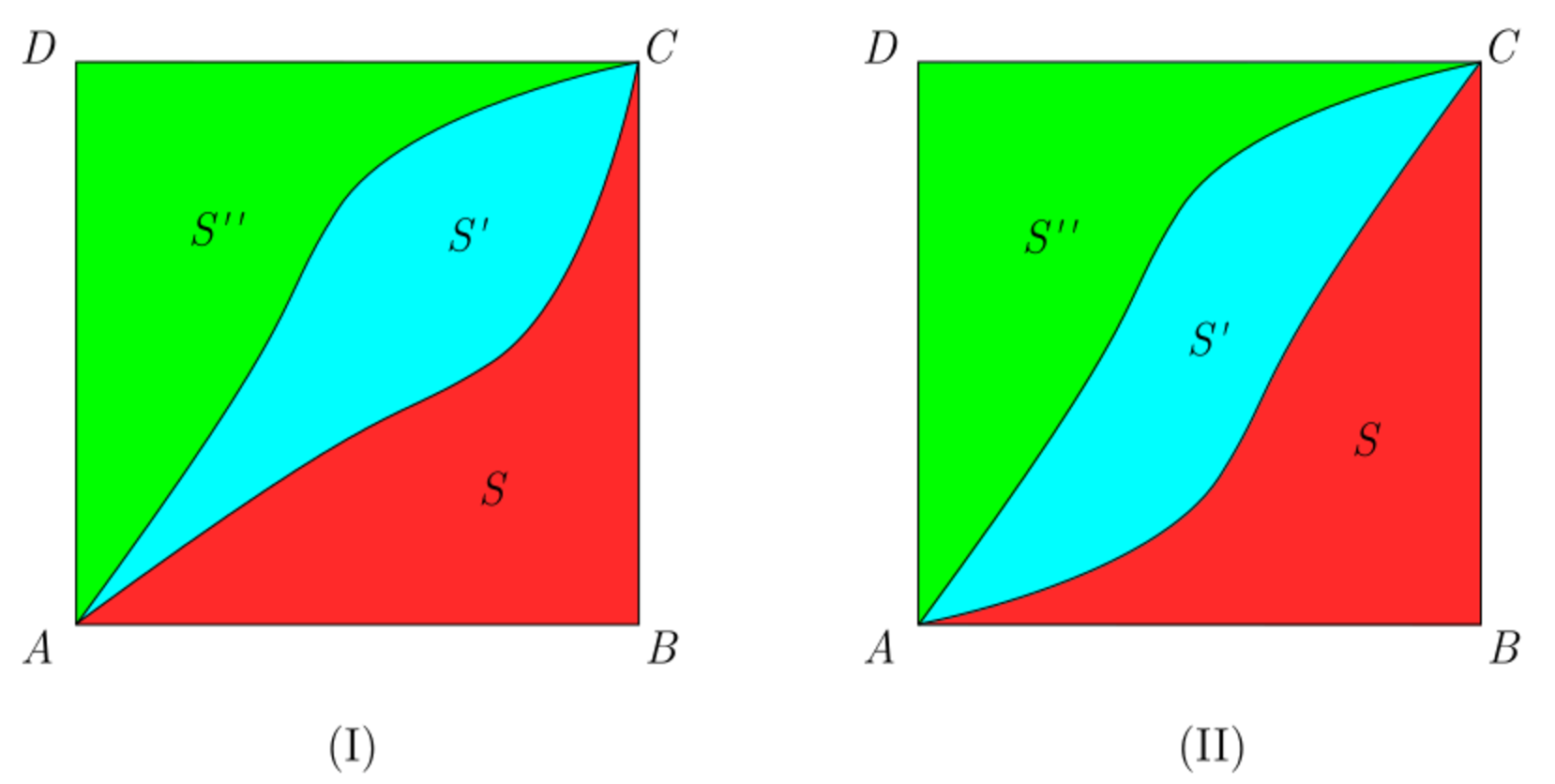}
\end{figure}

\begin{theorem}
There are no trisected squares with congruent regions with a region of type 2(e) (a region with exactly 2 insertion points such that the insertion points are on the same edge).
\end{theorem}
\begin{proof}
Suppose that there exists a trisected square with a region of type 2(e). Then it contains a region $S$ with exactly two insertion points, G and H, as shown in the figure below.  Consider the common boundary of the other two regions $S'$ and $S''$. We denote it as EF. Since neither $S'$ nor $S''$ can contain opposite corners of the rectangle, EF either connects \textbf{(I)} a point on the boundary of the rectangle with a point on $\partial S$ that is not an insertion point, or \textbf{(II)} a point on the boundary of the rectangle with another point on the boundary of the rectangle. 

\textbf{(I)}:
We will proceed by reductio ad absurdum. Consider the cases in which AB maps to a line segment in S. Note that the two segments AE and BH are bounding segments (if you extend them to lines the shape won't go past either segment). 
Thus, if AB maps to anywhere on curve GE, we would get that BH and AE map to CG and DE (in no particular order) because shape S spans the whole height of the rectangle and those two segments must bound it. This is clearly impossible because this would make shape S a rectangle.

Next, we will look the case where AB maps to DE. We know by definition that AB=DE, and by symmetry CD must map to EA, hence CD=EA. Note that both HB and CG must map to EF and only EF. Therefore, GF must map to HF and only HF, hence GF=HF. However, this implies that S' is symmetric, which in turn implies S and S" is symmetric. The only axis of symmetry that S can have is from D to the midpoint of GH, so GF=GH. This implies that S' is a triangle(a equilateral triangle to be exact), while S and S" are pentagons. Therefore, S and S" cannot be congruent to S', so AB cannot map to DE.

Now we will look at the case where AB maps to CD. We break this up into cases based on where in shape S segment AE maps. If segment AE maps to segment DE, then shape S'' is mapped to S by a reflection over the perpendicular bisector of DA. This would imply that the perpendicular bisector of DA is a bounding line for both S and S''. This is impossible by the same argument in the previous case. Next, if AE maps to CG, then we have HB maps to DE. However, this would imply HB=DE and CG=AE, meaning CG+HB=DE+AE=DA, and thus GH=0, which is impossible. Thus, in this case, we have the the perpendicular bisector of DA bounds both regions S and S''. Therefore, in order for EF to be the border between the two shapes, EF must be the perpendicular bisector of DA (or, rather, part of it). However, this forces region S' to have 0 area while S and S" have half the area of ABCD, which is clearly impossible.

Finally, we will look at the case where AB maps to CG. We must consider two cases: \textbf{(i)} EA maps to EG, or \textbf{(ii)} EA maps to CD.
\textbf{(i):} If EA maps to EG, then evidently EG is a straight line perpendicular to AD. Similarly, BF maps to DE, so BF is also a straight line. Therefore, S' again has 0 area while S and S" have half the area of ABCD. This is again impossible.
\begin{figure} [h!]
    \centering
	\includegraphics[width=0.6\textwidth]{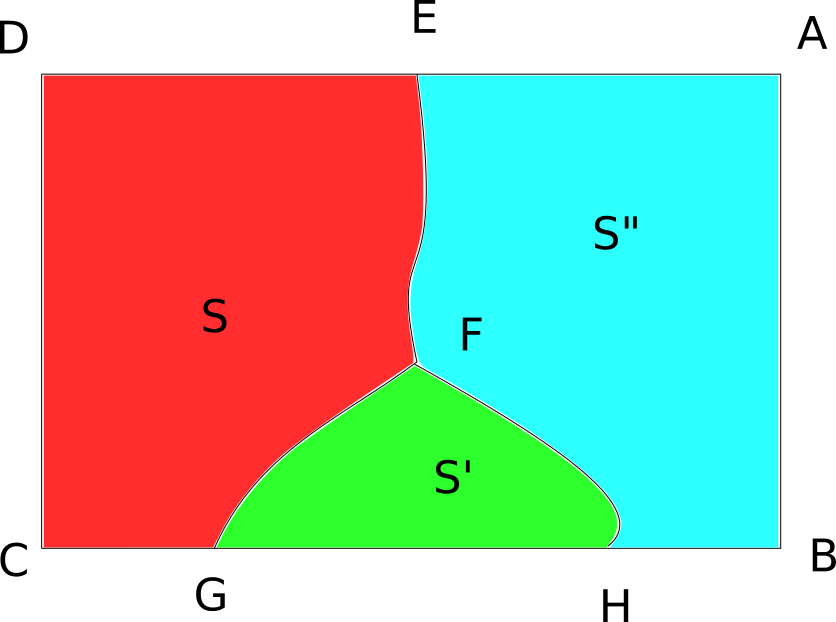}
\end{figure}
\newpage
\textbf{(II)}: There are three ways this is possible: \textbf{(i)} EF connects a point on $AD$ with a point on $CD$, \textbf{(ii)} EF connects a point on $AD$ with a point on $AG$ but not $G$, \textbf{(iii)} EF connects a point on $AD$ with $G$ (the other cases are symmetrical). Remember, $S'$ or $S''$ cannot contain opposite corners of the rectangle. 

\textbf{(i)}: $\partial S$ has $1+m$ straight line segments, $\partial S''$ has $3+n$ straight line segments, and $\partial S'$ has $4+m+n$ straight line segments, where $m$ and $n$ are non-negative integers. However, there are no non-negative integer solutions to $1+m=3+n=4+m+n$, so the three regions cannot be congruent. \textbf{(ii)} has the same argument. \textbf{(iii)} is also similar, except $\partial S'$ has $3+m+n$ straight line segments. 

It follows that a trisected square cannot have a region of type 2(e).
\end{proof}
\begin{theorem}
There are no trisected squares with congruent regions with a region of type 2(f) (a region with exactly 2 insertion points such that the insertion points are on a corner and an edge adjacent to that corner).
\end{theorem}
\begin{proof}
Suppose that there exists a trisected rectangle with a region of type 2(f). Then it contains a region $S$ with exactly two insertion points as shown in the figure below. Since $S$ does not contain opposite corners of the rectangle, the other two regions, $S'$ and $S''$, cannot contain opposite corners. So either (I) $A$ and $D$ belong to $S'$ and $B$ and $C$ belong to $S''$, or (II) $A$ and $B$ belong to $S'$ and $C$ and $D$ belong to $S''$. 

\textbf{(I)}: $S'$ must contain the entire edge $AD$ and $S''$ must contain the entire edge $BC$. $S'$ must contain a portion of $AD$ adjacent to $A$, $S''$ must contain a portion of $BC$ adjacent to $B$, and $S''$ must contain a portion of $CD$ adjacent to $C$. Furthermore, $S'$ must contain a portion of $\partial S$ adjacent to $D$, not on the boundary of the rectangle. There is an insertion point $P_1$ of $S'$ on $AB$. $S'$ either (i) has an insertion point on $CD$ or (ii) does not have an insertion point on $CD$. 

\textbf{(i)}: Suppose that $\partial S''$ contains $3+n$ line segments and $\partial S$ contains $1+m$ line segments. Then $\partial S'$ contains $3+m+n$ line segments. There are no non-negative integer solutions to $3+n=1+m=3+m+n$, so the three regions cannot be congruent. 

\textbf{(ii)}: We will start by noting a key attribute about figure $S''$. When looking at segment BC in figure, $S''$, the two straight segments that come off of it at right angles are bounding segments. By this, we mean that if we extended each of these segments on forever, region $S''$ would never cross through (go to the other side, not just intersect) either of the aforementioned segments. We will proceed by casework on where is region $S$ segment $BC$ maps. 
	First of all, we will look at whether $BC$ can map to the segment of region $S$ on the boundary of the rectangle. If this were the case, because because both segments coming out of $BC$ are straight bounding segments coming out at a right angle, it is clear that one of these two segments must map onto $AD$. However, this contradicts our construction of region $S$, thus, making this case impossible.
    Secondly, we will look at whether $BC$ can map to somewhere on the boundary of region $S$ not on the boundary of the trisected rectangle. If this were the case, because both segments coming out of $BC$ are straight and bounding segments, $P_1$ and the intersection of regions $S$, $S''$, and the trisected rectangle's boundary must map to the intersections of regions $S$ and $DC$. This would imply that shape S is made entirely of straight segments along the border,has 5 straight edges on the border,  and it clearly must be convex. However, it is clearly impossible for this to be true, shape S' to have only straight border segments, shape S' to have 5 straight segments on its border, and for shape S' to be convex.

\textbf{(II)}: $S'$ must contain the entire edge $AB$. Furthermore, it must contain a portion of $AD$ adjacent to $A$ and a portion of $BC$ adjacent to $B$. Then $S'$ has an insertion point $P_1$ on $AD$ and $P_2$ on $BC$. Consider the part $Z$ of $\partial S'$ connecting $P_1$ and $P_2$ not along the boundary of the rectangle. If $Z$ intersects $\partial S'$, there would be more than 3 regions, so $Z$ does not intersect $\partial S'$. Let $3+n$ be the number of line segments along $\partial S'$ and $1+m$ be the number of line segments along $\partial S$. Then $\partial S''$ has $3+m+n$ line segments. There exist no nonnegative $m,n$ such that $3+n=1+m=3+m+n$, so the three regions cannot be congruent.

\begin{centering}
\includegraphics[width=0.7\textwidth]{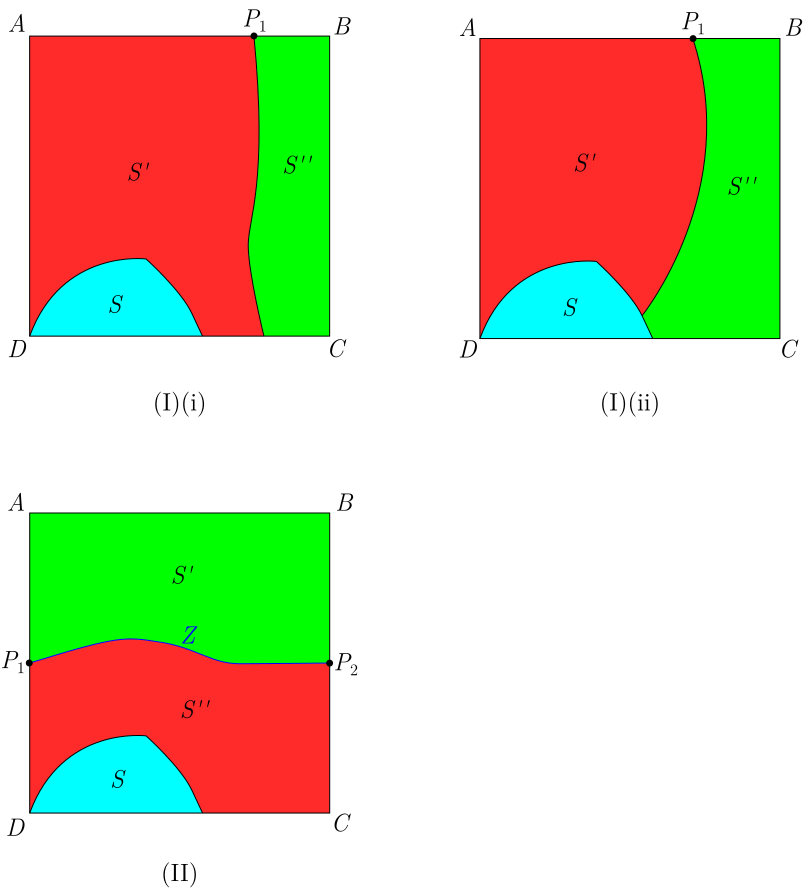}
\end{centering}

\end{proof}
\begin{theorem}
There are no trisected rectangle with congruent regions with a region of type 2(g) (a region with 2 insertion points such that the insertion points are on a corner and an edge opposite to that corner).
\end{theorem}
\begin{proof}
\begin{definition} Let $\sigma(XY)$ for curve $XY$ be the number of straight line segments not contained in any other line segment that are contained in curve $XY$. Similarly, let $\sigma(S)$ for shape $S$ be the number of straight line segments not contained in any other line segment contained in the perimeter of shape $S$.
\end{definition}
We start by observing that once we know that we have a region of type $2(g)$, any on curve will uniquely determine the remaining two regions and vice versa. We will all this curve the partition. To reduce the number of cases, we can also note that if any region contains opposite corners, than it contains a straight line segment inside of it that is bigger than any line segment contained in the region of type $2(g)$. Thus, we consider only cases in which no region contains opposite corners. Let the type $2(g)$ region in each diagram be called region $S$. Also, note that whenever we say curve $AB$ can't map to curve $CD$, we also mean that curve $AB$ can't map to any curve contained in $CD$. Thus, we have the following cases:

$a)$ We let the partition go from the insertion point of region $S$ not on a corner of the trisected rectangle to a point on the opposite edge and not on a corner, as in the diagram below.\\

\begin{centering}
\includegraphics[width=0.5\textwidth]{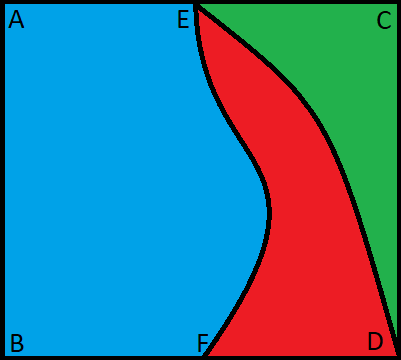}\\
\end{centering}

Note that $\sigma(ABFE)=\sigma(FED)=\sigma(CED)$, $\sigma(ABFE)=3+\sigma(FE)$, $\sigma(FED)=1+\sigma(FE)+\sigma(DE)$, and that $\sigma(DEC)=2+\sigma(DE)$. Thus, $3+\sigma(FE)=1+\sigma(FE)+\sigma(DE)\implies \sigma(DE)=2$ and $1+\sigma(FE)+\sigma(DE)=2+\sigma(DE)\implies \sigma(FE)=1$. Observe that $AB$ can't map to $EF$ without forming a rectangle because $\sigma(FE)=1$. Next, because $\measuredangle BDE < 90^{\circ}$, $AB$ can't map to $FD$. Finally, observe that because $\sigma(DE)=2$ and $0^{\circ} < \measuredangle BDE < 90^{\circ}$, curve $AB$ can't map to curve $DE$. Thus, case $a$ is impossible.

$b)$ We let the partition go from the insertion point of region $S$ not on a corner of the trisected rectangle to a point on the opposite corner, as in the diagram below.

\begin{centering}
\includegraphics[width=0.8\textwidth]{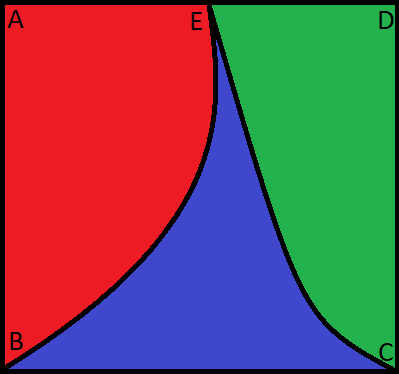}\\
\end{centering}

Observe that $\sigma(EB)+2=\sigma(CE)+2=\sigma(BE)+\sigma(CE)+1 \implies \sigma(BE)=\sigma(CE)=1$. However, note that region $AEB$ has a right angle, so region $BEC$ must as well. The only way for this to happen is if two straight segments contained in $BE$ and $EC$ meet at a right angle at point $E$. However, if this is the case note that region $AEB$ has all three of its straight segments connected.If this were true, then we would see that $BE$ must be just one straight line segment as must be $EC$. However, this can't be true because then we would get $BC>BE$ because $BC$ is the hypotenuse of $BEC$ and $BE>BA,AE$ because $BE$ is the hypotenuse of $BAE$. Thus, our regions aren't congruent because the Boundary of $BEC$ contains a line segment longer than any in the border of $BAE$.

$c)$ We let the partition go from the insertion point of region $S$ not on a corner of the trisected rectangle to a point on the adjacent edge and not on a corner, as in the diagram below.\\

\begin{centering}
\includegraphics[width=0.8\textwidth]{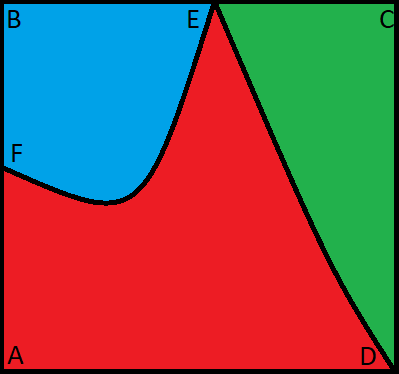}\\
\end{centering}

Observe that $\sigma(FE)+2=\sigma(ED)+2=\sigma(FE)+\sigma(DE)+2 \implies \sigma(FE)=\sigma(DE)=0$.Thus, we see that right angle $\angle FAD$ maps to right angle $\angle EBF$ and  to $\angle ECD$. Thus, the length of the non straight curves on the boundary of each shape must be equal so $FE=ED=FE+ED\implies FE=ED=0$, which is clearly impossible.

$d)$ We let the partition connect opposite edges of the trisected rectangle without ending on a corner of the trisected rectangle or an insertion point of our type $2(g)$ region, as in the diagram below.\\

\begin{centering}
\includegraphics[width=0.8\textwidth]{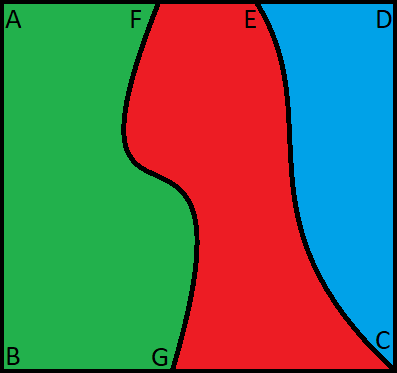}\\
\end{centering}

Observe that $\sigma(FG)+3=\sigma(EC)+2=2+\sigma(EC)+\sigma(FG)\implies \sigma(FG)=0, \sigma(EC)=1$. However, in that, case region $EFGC$ clearly has fewer right angles than region $AFGB$.\\
$e)$ We let the partition go from the corner of the trisected rectangle adjacent to the insertion point of the $2(g)$ on a corner of the trisected rectangle to the edge opposite that point not on an insertion point of the type 2(g) region or on a corner of the trisected rectangle, as seen in the diagram below.\\

\begin{centering}
\includegraphics[width=0.8\textwidth]{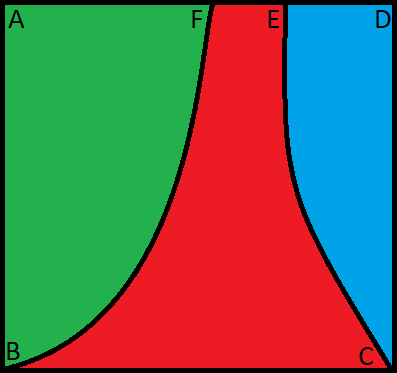}\\
\end{centering}

Observe that $\sigma(FB)+2=\sigma(FB)+\sigma(ED)+2=f+(EC)\implies \sigma(EC)=\sigma(FB)=0$. However, then region $FECB$ has no right angles while region $AFB$ has one right angle. \\
$f)$ We let the partition go from the edge opposite the 2(g) (the one it doesn't touch at all) but not on the endpoint of that edge to the edge of the trisected rectangle adjacent to it but not on that edge's endpoints, as seen in the diagram below.\\

\begin{centering}
\includegraphics[width=0.8\textwidth]{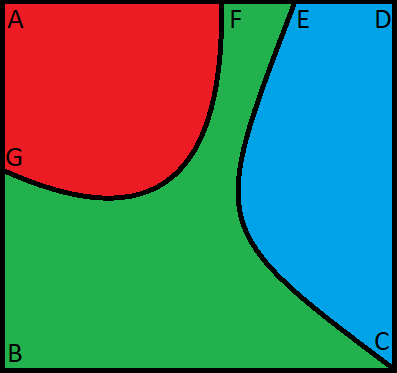}\\
\end{centering}

Observe that $\sigma(FG)+2=\sigma(EC)+2=\sigma(EC)+\sigma(FG)+3\implies \sigma(FG)=\sigma(EC)=-1$, however, $\sigma(XY)$ for any segment $XY$ clearly must be non-negative.\\
$g)$We let the partition go from the boundary of the trisected rectangle on the same side that the type $2(g)$ region has an insertion point on not on the insertion point of the type $2(g)$ region or a corner of  the square to the insertion point of the type $2(g)$ region that is on the boundary of the trisected rectangle, as seen in the diagram below.\\

\begin{centering}
\includegraphics[width=0.8\textwidth]{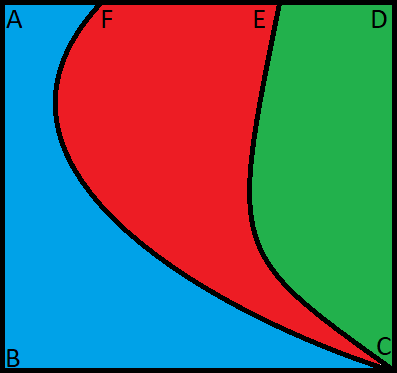}\\
\end{centering}

Observe that $\sigma(FC)+3=\sigma(FC)+\sigma(CE)+1=\sigma(CE)+2\implies \sigma(CE)=2,\sigma(FC)=1$. We will proceed by casework on what part of the boundary of region $DEC$ segment $AB$ maps to. First, $AB$ clearly can't map to $DC$ because we clearly can't have a right angle at point $C$. Next, if $AB$ maps to $DE$ , then $AD>AB$. However, This would mean that the segment that $BC$ maps to must be perpendicular to $AD$ , therefore parallel to $DC$, but it mus also be longer than $AB$, which is clearly impossible. Finally, if $AB$ maps to some segment in curve $EC$,then we have  that $AF$ maps to $DE$ meaning that $AB$ maps to a segment parallel to $DC$ with the same length as $DC$, which is impossible with the given conditions and would cause there to be a rectangle.\\

$h)$We let the partition go from the boundary of region $S$ not on one of its insertion points to the opposite edge of the trisected rectangle, as seen  in the diagram below.\\

\begin{centering}
\includegraphics[width=0.99\textwidth]{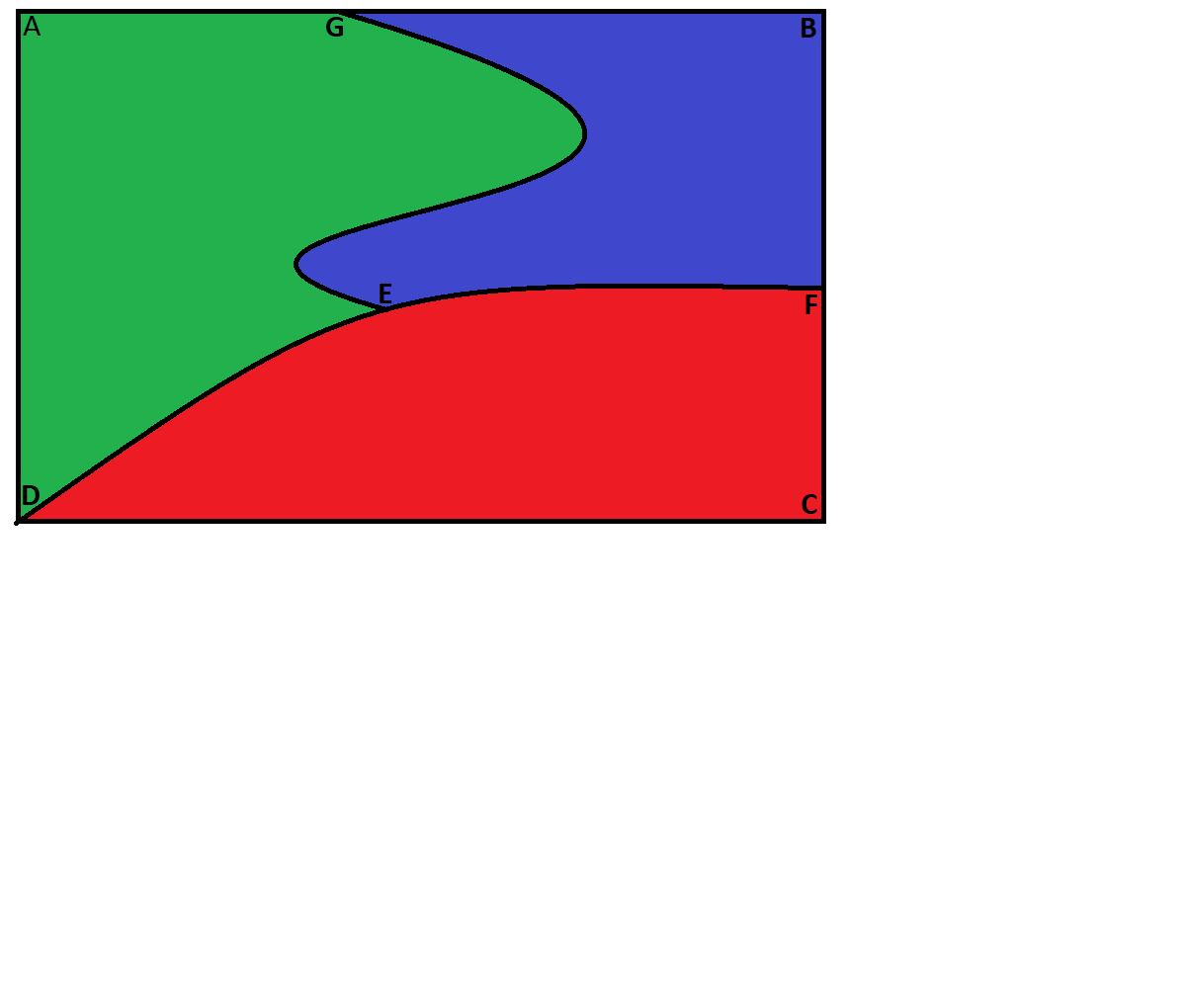}\\
\end{centering}

Assume, without loss of generality, that $CD\geq AD$ (if this is not the case, merely rotate $ABCD$ $90 \degree$ counter-clockwise about point $D$). We intend to prove by contradiction that $\angle FCD$ can't map to $\angle GAD$. Assume for the sake of contradiction that$\angle FCD$ maps to $\angle GAD$. Note that $CD>AG$, so $CD$ maps to $AD$, and $ABCD$ is a square.  However, then $CD$ can't map to $BG$  or $BF$ because both are shorter. Because both segments from $C$ at a right angle going out to $D$ and $F$ are connected, $D$ mus map to $G$ and $F$ must map to $F$ or vice versa . Thus, each region consists of four straight lines, two opposite pairs of line segments which form right angles.  Thus, we can observe that  each shape must form a right angle at F, so they sum up to $270 \degree$, but they should sum $360 \degree$, so this is a contradiction.  
Thus, these two angles can't map to each other, and the right angle $\angle FCD$ must map to $F$, so the center angles add up to $270 \degree $ not $360 \degree$ as described in the previous paragraphs, meaning this case is impossible.

$i)$ We let the partition go from the boundary of region $S$ not on one of its insertion points to the corner that is adjacent to the insertion point of region $S$ in one of the corners of the rectangle, as seen in the diagram below.

\begin{centering}
\includegraphics[width=0.8\textwidth]{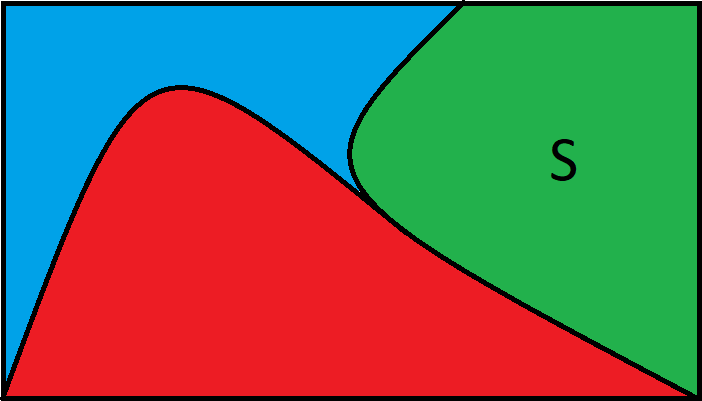}\\
\end{centering}

This case is  impossible because there is a region of type $2(b)$, which was proved impossible in theorem $7$.

\end{proof}
\begin{corollary}
There are no trisected squares with congruent regions of type 2.
\end{corollary}
\begin{proof}
Follows theorems 6-12.
\end{proof}
\begin{corollary}
A trisected square is impossible.
\end{corollary}
\begin{proof}
Follows corollaries 1-3.
\end{proof}
\section{Conclusion}

In this paper, we reviewed the history of the paper, and come up with a alternate proof that shows that this problem is not solvable. However, this is not the end to this string of problems. This is a part of a larger problem, can a rectangle be divided into three congruent, non-rectangular sides, as a square is classified as a rectangle. We plan to edit our proof in the future to address this problem. Beyond this lies a another set of problems waiting to have their answers unveiled. Some of these may much more complex and intriguing than this problem, such as trying to divide a rectangle into nine parts, or trying to trisect a  11-sided figure, however these would most likely need to be proved through a computer, due to the sheer amount of possibilities.

\bibliographystyle{abbrv}

\end{document}